\documentclass[11pt]{article}

\usepackage{amssymb}
\usepackage{amsfonts, eufrak}
\usepackage[dvips]{graphics}
\usepackage[latin1]{inputenc}
\usepackage{amsmath}

\def\dd{\displaystyle}

\def\ca{\c{c}\~ ao }

\def\la{\langle}
\def\ra{\rangle}

\def\dd{\displaystyle}
\def\eps{\epsilon}

\def\tto{\longrightarrow}

\def\rr{\mathbb{R}}

\newtheorem{thm}{Theorem}

\begin{document}
\title{A Note on Maximal Averages in the Plane }
\author
    {Jos\' e A. Barrionuevo%
     \thanks{josea@mat.ufrgs.br} \\
     Lucas Oliveira%
     \thanks{lucas\_ gnomo@hotmail.com} \\
     Departamento de Matem\' atica \\
     UFRGS \\
     Av. Bento Gon\c calves 9500, 91509-900  Porto Alegre, RS, Brasil.}

\date{}
\maketitle

\section{Introduction}
The study of directional maximal operators goes back at least to the works of C\' ordoba \cite{C1} and
Str\" omberg \cite{Str1} in the 70's and is as recent as the results of Alfonseca, Soria and Vargas
\cite{AVS1}, \cite{AVS2},
Karagulyan and Lacey \cite{KL} and Bateman and Katz \cite{Ba1}, \cite{BK1}. Progress in the area involved
the work of many authors and we refer to \cite{W}, \cite{AVS2}, \cite{KL} and references therein
for historical background. \\
The general situation considered in \cite{AVS2} is the following. Suppose $\dd \Omega = \Omega_0\cup_j
\Omega_j \subset [0,1]$ is a closed set of slopes in the plane, where
$\Omega_0 = \{ \theta_1 > \theta_2 > \cdots > \theta_j > \cdots \}$ is an ordered subset of $\Omega$ and
$\Omega_j = [\theta_{j-1}, \theta_j) \cap \Omega$. Denote by $\dd \mathcal{B}_{\Omega}$ be the collection of
all rectangles so that its longest side has slope in $\Omega$, and define
\[ M_{\Omega} f(x) = \sup_{x\in R\in \mathcal{B}_{\Omega}} \frac{1}{|R|} \int_R|f(y)|\, dy \]
with similar definitions for $M_{\Omega_0}, M_{\Omega_i}$. If $||\, L\, || = \sup\{\, ||\, L f\, ||:\, || f ||
\leq 1\}$ represents the norm of the linear or sub-linear operator $L$ acting in $L^2(\rr^2)$, the following
is proved in \cite{AVS2}
\begin{thm}\label{Thm1}
The following inequality holds
\begin{equation} \label{AVS}
 ||M_{\Omega}|| \leq \sup_i \left|\left| M_{\Omega_i}\right|\right| + C||M_{\Omega_0}||
 \end{equation}
 where $C$ is an absolute constant.
 \end{thm}
 It is known that this inequality implies the $\log N$ upper bound for the norm of the operator associated
 with $N$ arbitrary directions first proved in \cite{K1}, as well as the case of $N$-lacunary directions,
 \cite{Ba1}, \cite{KL}. \\
 In the next section we provide a self-contained proof of (\ref{AVS}) by the $TT^{\ast}$ method in the
 lines of \cite{B1}, \cite{B2}. General description of the method can be found in \cite{St},
 \cite{T1} and \cite{B1}. In a certain sense these $TT^{\ast}$ arguments can be seen as another instance
 of the Tensor Power Trick as presented in \cite{T2}. Note that the case of $\Omega_0$ lacunary and each
 $\Omega_j$ uniformly distributed in the sector was consider in \cite{C2} and \cite{B1} where sharp weak
 and strong $L^2$ estimates were obtained.

 As a second illustration of the method we prove the following slight improvement of C\' ordoba´s original
 estimate \cite{C1}. Let $M_{\delta}$ denote the maximal operator with respect to the basis of $h\times \delta h$
 rectangles. It is known that $||\, M_{\delta}\, || \leq C |\log \delta |$. Define the following {\em grand}
 maximal operator
 \begin{equation}\label{GM}
 \mathcal{GM} f(x) = \sup_{0 < \delta < 1/2}\, \frac{1}{|\log\delta |}\, {M}_{\delta}\, f(x)
 \end{equation}
Then we prove
\begin{thm}\label{thm2}
 $\mathcal{GM}$ is bounded on $L^2(\rr^2)$.
 \end{thm}
 Observe that the factor $\frac{1}{|\log\delta\, |}$ in (\ref{GM}) is sharp. For if it is replaced by
 $A_{\delta}$ satisfying \\
 $\limsup_{\delta\to 0^+} A_{\delta}\, |\log\delta | = \infty$, then $\mathcal{GM}$ is unbounded. \\
 In what follows $C$ denotes a constant not necessarily the same on each occurrence but which is absolute.
 Norms of vectors $|| f ||, ||T f ||$ are Lebesgue spaces $L^2(\rr^2)$ norms while norms $|| S ||$ of
 operators are the uniform operator norm.

\section{Proof of Theorem 1}
Given a measurable map $\dd \Phi: x\in\rr^2 \tto R_x\in \mathcal{B}_{\Omega}$ define the linear operator $T$ associated to $\Phi$ by
\begin{equation}\label{T}
 Tf(x) = \frac{1}{|R_x|}\int_{R_x} f(y) \, dy = \int_{\rr^2} k(x,y)\, f(y)\, dy
\end{equation}
where the kernel is $\dd k(x,y) = \frac{1}{|R_x|}\, \chi_{R_x}(y)$.

We define the following auxiliary operators: if $R_x$ has slope in $\Omega_j$, let $R^{\prime}_x,\,
R^{\prime\prime}_x$ be rectangles with same center and twice the dimensions of $R_x$ but with its longest side with slopes $\theta_{j-1}$ and $\theta_j$ respectively. Define
\begin{equation}\label{T0}
T_0 f(x) = \frac{1}{|R_x^{\prime}|}\int_{R_x^{\prime}} f(y) \, dy + \frac{1}{|R_x^{\prime\prime}|}\int_{R_x^{\prime\prime}} f(y) \, dy
\end{equation}
and if $\widetilde{R}_x$ is the rectangle with same center and slope and  twice the dimensions of $R_x$, define $\widetilde{T}$ similarly to (\ref{T}). Observe that
\begin{equation}
|T_0f(x)| \leq T_0|f|(x) \leq C M_{\Omega_0} f(x)
\end{equation}
\begin{equation}\label{T1}  |Tf(x)| \leq T|f|(x) \leq M_{\Omega} f(x) \end{equation}
\begin{equation}\label{tT}  |\widetilde{T} f(x)| \leq \widetilde{T} |f|(x) \leq M_{\Omega} f(x) \end{equation}
and that for $f\geq 0$,  $Tf(x) \leq C\, \widetilde{T}f(x)$. The adjoint of $T$ is given by
\[ T^{\ast}f(y) = \int_{\rr^2} \frac{1}{|R_x|}\, \chi_{R_x}(y) f(x)\, dx,  \]
and thus a straightforward calculation gives
\[ TT^{\ast} f(x) = \int_{\rr^2} \mathbb{K}(x,z)\, f(z)\, dz \] where
\begin{equation}\label{KK}
 \mathbb{K}(x,z) = \int_{\rr^2} \frac{1}{|R_x|}\, \chi_{R_x}(y)\, \frac{1}{|R_z|}\, \chi_{R_z}(y)\, dy
 = \frac{|R_x\cap R_z|}{|R_x| |R_z|} \end{equation}

Define $A_j = \{ x\in\rr^2 : R_x \in \mathcal{B}_{\Omega_j}\}$. Then the $A_j$ are measurable, disjoint and cover $\rr^2$. We write
\[ \mathbb{K}(x,z) = \sum_{i,j} \chi_{A_i}(x)\, \mathbb{K}(x,z)\, \chi_{A_j}(z) = \mathbb{K}_1(x,z) +
\mathbb{K}_2(x,z) \]
where $\dd \mathbb{K}_1(x,z) =  \sum_{i}  \chi_{A_i}(x)\, \mathbb{K}(x,z)\, \chi_{A_i}(z)$, and for $s=1, 2$ define the associated operators
\[ \mathbb{T}_s f(x) = \int_{\rr^2} \mathbb{K}_s(x,z) f(z)\, dz \]
Since $ ||\, \chi_{A_i}\, T\, || = ||\,(\chi_{A_i}\, T)^{\ast}\, || \leq ||\, M_{\Omega_i}\, ||$, the disjointness of $A_{i's}$ yields
\begin{eqnarray*}
||\, \mathbb{T}_1 f\, ||^2 & = & \sum_i\, ||\, \chi_{A_i}\, (\chi_{A_i} T)\, (\chi_{A_i} T)^{\ast}\, \chi_{A_i} f\, ||^2 \\
\nonumber  & \leq & \sum_i \left|\left| M_{\Omega_i}\right|\right|^4\, \left|\left| \chi_{A_i} f \right|\right|^2 \\
 & \leq & \sup_i \left|\left| M_{\Omega_i}\right|\right|^4\,\,\sum_i \left|\left| \chi_{A_i} f \right|\right|^2 \\
 & = & \sup_i \left|\left| M_{\Omega_i}\right|\right|^4\,\, ||f||^2
\end{eqnarray*}
establishing that
\begin{equation}\label{T1}
||\, \mathbb{T}_1\, || \leq \sup_i\, ||\, M_{\Omega_i}\, ||^2
\end{equation}
To estimate $\mathbb{T}_2f$, we dominate (\ref{KK}) by replacing one of the rectangles by one having slope in $\mathcal{B}_0$, that is, there are $R_x^{\prime}, R^{\prime}_z\in \mathcal{B}_0$ such that
\begin{equation}\label{geom}
\frac{|R_x\cap R_z|}{|R_x| |R_z|} \leq C\, \max \left\{ \frac{|\widetilde{R_x}\cap R^{\prime}_z|}{|\widetilde{R_x}| |R^{\prime}_z|}
,\,  \frac{|R^{\prime}_x\cap \widetilde{R_z}|}{|R^{\prime}_x| |\widetilde{R_z}|}\right\}
\end{equation}
This implies that for all $x$ in $\rr^2$ and $f\geq 0$
\begin{eqnarray*}
 \mathbb{T}_2 f(x)  & \leq & C\, ( \widetilde{T} T_0^{\ast} f(x) + T_0\widetilde{T}^{\ast} f(x) )
\end{eqnarray*}
therefore we obtain the following pointwise inequality for $TT^{\ast}$
\begin{equation}\label{TT}
TT^{\ast} f (x) \leq \mathbb{T}_1 f(x) + C\, ( \widetilde{T} T_0^{\ast} f(x) + T_0\widetilde{T}^{\ast} f(x) )
\end{equation}
valid for non negative functions. By (\ref{T1}) and (\ref{TT}) we have the following operator´s norm inequality
\begin{equation}\label{TT2}
||\, TT^{\ast}|| \leq \sup_i \left|\left|\, M_{\Omega_i}\right|\right|^2 + C\, ||\, \widetilde{T}||\,\, ||\, T_0||
\leq \sup_i \left|\left|\, M_{\Omega_i}\right|\right|^2 + C\, ||\, M_{\Omega}||\,\, ||\, M_{\Omega_0}||
\end{equation}
Now given $\eps > 0$ first select $f\geq 0$ in $L^2$ with $|| f || = 1$ such that
\[ (1+\eps )\, ||\, M_{\Omega} f || \geq ||\, M_{\Omega}\, ||. \]
Next, choose $R_x$ in (\ref{T}) satisfying
\[ (1+\eps )\, Tf(x) \geq M_{\Omega} f(x) \]
and observe that this implies
\begin{equation}\label{Rx}
(1+\eps )^2\, ||\, T\, || \geq ||\, M_{\Omega}\, ||.
\end{equation}
Finally pick $g\geq 0$ in $L^2$ with $|| g || = 1$ so that
\begin{equation}\label{g}
(1+\eps )\, ||\, T^{\ast} g || \geq ||\, T^{\ast} || = ||\, T\, ||. \end{equation}
By (\ref{Rx}), (\ref{g}), and Cauchy-Schwarz, we get
\begin{equation}\label{MM}
||\, M_{\Omega}\, ||^2  \leq  (1+\eps )^6\, ||\, T^{\ast} g\, ||^2
  =  (1+\eps )^6\,  \la\, TT^{\ast} g, g\,\ra
  \leq  (1+\eps )^6\, ||\,  TT^{\ast}\, ||
 \end{equation}
Letting $\eps \to 0$ in (\ref{MM}) and inserting the result into (\ref{TT2}) immediately yields (\ref{AVS}) $\blacksquare $
\section{Proof of Theorem 2}
We use a variant of Proposition 1 of \cite{B1}. First note that by the monotone convergence theorem it suffices to prove
\begin{equation}\label{M0}
||\, \mathcal{GM}_0\, || = || \,  \sup_{0<\delta_0 < \delta < 1/2}\, \frac{1}{|\log\delta |}\, {M}_{\delta}\, || \leq C
\end{equation}
with $C$ independent of $\delta_0$. We will also assume that the set of directions is restricted to
$\theta\in [0,1/10]$, thus making the $y$-vertical direction uniformly transversal to all directions and with the cost of a multiplicative constant in the operator norm.
For each vector $m = (h,\theta , \delta )$, let $R_m$ be the $h\times \delta h$ rectangle centered in $x$ with longest side pointing along $\theta$. Let $l_m = (1+| \log\delta |)$ and define the positive self-adjoint linear operator $T_m$ by
\begin{equation}
T_m\, f(x) =  l_m^{-1}\,\frac{1}{|R_m|}\int_{R_m} f(y)\, dy \end{equation}
If $m = (h,\theta , \delta ),\, n = (k,\beta , \eta )$ with $h > k$, a direct calculation (carried out in \cite{B1}) gives
\begin{equation}\label{TT3}
T_mT_n f(x) \leq C\, ( l_m\, l_n)^{-1}\, \frac{1}{|R|}\,\int_R f(y)\, dy
\end{equation}
where $R$ is a rectangle of dimensions $2h\times 2w$ with  $w =  \max \{ \eta\,k,  \delta\, h,  k\, \sin |\theta - \beta |\}$ with longest side ($2h$) pointing along $\theta$. We consider two cases separately. \\
{\bf case 1:} $w = \delta\, h$. \\
Define the linear positive self adjoint operator $S_m$ by
\[ S_m f(x) = l_m^{-1}\,\frac{1}{|\widetilde{R}_m|}\int_{\widetilde{R}_m} f(y)\, dy \]
where $\widetilde{R}_m$ is a rectangle with same center and slope of $R_m$ and twice the dimensions. Let
$S f(x) = \sup_m |S_m f(x)|$ be the corresponding maximal operator. Note that since for $f\geq 0$, \\
$ S_m f(x) \leq  \mathcal{GM}_0 f(x)  $, we have
\begin{equation}\label{caso1}
 ||\, S\, || \leq ||\, \mathcal{GM}_0 \, ||
\end{equation}
and in this case (\ref{TT3}) becomes
\begin{equation}\label{TTcaso1}
T_mT_n f(x) \leq C\, ( l_n)^{-1}\, S_m f(x) \leq C\,  S_m f(x)
\end{equation}
{\bf case 2:} $w = \max \{ \eta\,k,  k\, \sin |\theta - \beta |\}$ \\
We have $w\in (\eta\, k\, 2^{j-1},  \eta\, k\, 2^j]$ for some $j\leq 2\, l_n$. Define
\[ W_n\, f = (2\, l_n)^{-1}\, \sum_{j=1}^{2l_n} H_{n, j}\; f\]
where for each $n = (k, \beta , \eta )$ and $j \leq 2\, l_n$,
\[ H_{n, j}\, f(x) = (\eta\, k\, 2^{j+1})^{-1}\,\int_{-\eta\, k\, 2^j}^{\eta\, k\, 2^j}\, f(x + t(0,1))\,
 dt \]
Observe that each $H_{n, j}$ is positive, self adjoint and $|\, H_{n, j}\, f(x) \,| \leq C\, M^y f(x)$, where $M^y$ is the one dimensional
Hardy-Littlewood maximal operator acting in the $y$ direction. Thus for the maximal operator
$W\, f(x) = \sup_n|W_n f(x)|$, we have
\begin{equation}\label{caso2}
||\, W\, || \leq C\, ||\, M^y\, || \leq C
\end{equation}
yielding the following pointwise estimate for (\ref{TT3})
\begin{equation}\label{TTcaso2}
T_mT_n f(x) \leq C\; S_m\, W_n f(x) = C\; W_n\, S_m f(x) 
\end{equation}
Taking (\ref{TTcaso1}), (\ref{TTcaso2}) into account and considering also the case where $h\leq k$, we obtain
the final pointwise estimate for (\ref{TT3})
\begin{equation}\label{GG}
T_mT_n f(x) \leq C\,\,\left( S_m + S_n + S_m W_n + W_m S_n\right)\, f(x)
\end{equation}
For the reader´s convenience we adapt the argument in \cite{B1} to the present situation. Since all operators
involved are positive we only need to consider $f\geq 0$. Given $\eps > 0$ choose $f\in C_0\cap L^2(\rr^2)$
positive, with $||\, f\, || \leq 1$, then a measurable $m(x)$ such that
\begin{eqnarray}
(1+\eps )\, ||\, \mathcal{GM}_0\, f\, || & \geq & ||\, \mathcal{GM}_0\, ||  \\
(1+\eps )\, ||\, T_{m(\cdot )}\, || & \geq & ||\, \mathcal{GM}_0\, ||
\end{eqnarray}
Next choose $g(x)$ with $|| g || = 1$ satisfying
\begin{equation}\label{g}
 (1+\eps)\, ||\,(T_{m(\cdot )})^{\ast}\, g\, || \geq ||\,(T_{m(\cdot )})^{\ast}\, || = ||\,(T_{m(\cdot )})\, ||
\end{equation}
Note that, by positivity, inequality (\ref{GG}) reamins true if we take $m=m(x)$ to be any
measurable function. 
Plugging $g$ and $m(x)$ into (\ref{GG}) and using that positivity is preserved by taking adjoints we get
\begin{eqnarray}
T_n\, (T_{m(\cdot )})^{\ast}\, g(x) & \leq & C\,\,\left( (S_{m(\cdot )})^{\ast} + S_n + W_n \, (S_{m(\cdot )})^{\ast}
 + S_n\, (W_{m(\cdot )})^{\ast} \right)\, g(x) \\
  & \leq &  C\,\,\left( (S_{m(\cdot )})^{\ast} + S + W \, (S_{m(\cdot )})^{\ast}
 + S\, (W_{m(\cdot )})^{\ast} \right)\, g(x)
\end{eqnarray}
Observe that by (\ref{caso1}), (\ref{caso2}) we have the estimates
\begin{eqnarray}
||\, (S_{m(\cdot)})^{\ast}\, || = ||\, S_{m(\cdot )}\, || & \leq & ||\, \mathcal{GM}_0 \, || \\
||\, (W_{m(\cdot)})^{\ast}\, || = ||\, W_{m(\cdot )}\, || &  \leq & C\, ||\, M^y\, || \leq C
\end{eqnarray}
and that we can let $n = m(x)$ in (28) without changing the right hand side.
Taking inner product of (28) with $g$ and using Cauchy-Schwarz, (19), (21), (29), (30) we obtain
\begin{equation}\label{GG3}
||\,(T_{m(\cdot )})^{\ast}\, g\, ||^2 =  |\, \langle T_{m(\cdot )}\, (T_{m(\cdot )})^{\ast}\, g,\, g \rangle\, | \leq C \left(\, ||\, \mathcal{GM}_0 \, || +
||\, \mathcal{GM}_0 \, ||\; ||\, M^y\, ||\right) \leq C\, ||\, \mathcal{GM}_0 \, ||
\end{equation}
Taking (25) and (26) into (31) and letting $\eps\to 0$ we get the result.  $\blacksquare$

This theorem remains true if $M_{\delta}$ is replaced by the codimension 1 Kakeya operator, that is, the
maximal averages over $h^{n-1}\times \delta\, h$ rectangles acting on $L^2(\rr^n)$.

\end{document}